# Summation arithmetic functions with bounded terms, having a limit normal distribution law

## VICTOR VOLFSON

ABSTRACT The paper considers the properties of pseudo stationarity in a broad sense and pseudo strong mixing for sequences of random variables corresponding to arithmetic functions. Assertions on this topic have been proven. The implementation of these properties for known arithmetic functions has been verified. The article proves a statement about sufficient conditions under which a summation arithmetic function with bounded terms has a limit normal distribution law. The fulfillment of the specified sufficient conditions for known arithmetic functions is considered.

## 1. INTRODUCTION

An arithmetic function (in the general case) is a function $f : N \to C$ defined on the set of natural numbers and taking values on the set of complex numbers. The name arithmetic function is related to the fact that this function expresses some arithmetic property of the natural series.

A summation arithmetic function is a function:

$$S(x) = \sum_{n \le x} f(n). \tag{1.1}$$

The example of summation functions (1.1) is Mertens function:

$$M(\mathrm{x}) = \sum_{k \le x} \mu(k), \tag{1.2}$$

where $\mu(k)$ is Möbius function.

A Mobius function $\mu(\mathrm{k}) = 1$ if the natural number k has an even number of prime divisors of the first degree, $\mu(\mathrm{k}) = -1$ if the natural number k has an odd number of prime divisors of the first degree and $\mu(\mathrm{k}) = 0$ if the natural number k has prime divisors of not only first degree.

The study of summation arithmetic functions is usually carried out in two directions: the study of their distribution depending on the values of the natural argument and the study of their asymptotic behavior with the value of the natural argument $n \to \infty$.

The study of the distribution of arithmetic functions is usually reduced to the determination of their average values [1].

The probabilistic approach to the study of the distribution of arithmetic functions [2] provides new opportunities.

Probabilistic spaces are determined $\left(\Omega_n, \mathcal{A}_n, \mathbb{P}_n\right)$ by taking $\Omega_n = \{1, 2, ..., n\}$, $\mathcal{A}_n$ - all subsets $\Omega_n$, $P_n(A) = \dfrac{1}{n}\{N(m \in A)$, where $N(m \in A)$ is the number of members of the natural series that satisfy the condition $m \in A$ in this case.

Then an arbitrary (real) arithmetic function $f$ (or rather, its restriction on $\Omega_n$) can be viewed as a sequence of random variables $f_n$ defined on different probability spaces: $f_n : f_n(k) = f(\mathrm{k}), (\mathrm{k} = 1, ..., \mathrm{n})$.

Therefore, it is possible to determine the variance, other moments of higher orders, the distribution function and the characteristic function, sequences of random, corresponding to the arithmetic function.

If certain conditions are fulfilled, the specified sequence of random variables when the value $n \to \infty$ weakly converges in distribution to a certain distribution function. In this case, we say that the limit distribution of an arithmetic function is a given distribution function.

The theorem is proved [3] that the limit distribution of the arithmetic function - the number of prime divisors of a natural number when the value $n \to \infty$ is the function of the standard normal distribution. Conditions were found in [2], [4], [5] under which the limit distribution of additive and multiplicative functions have the standard normal distribution.

Definition. We understand as asymptotic independence of arithmetic functions that when the limit of the difference between the average value of the product of the arithmetic function $f$ and the product of the average values of the same function for different values of the argument when the value $n \to \infty$ tends to zero.

It was proved [6] that the property of asymptotic independence of the terms holds for the summation functions of Mertens $M(x) = \sum_{n \leq x} \mu(\mathrm{n})$ and Liouville $L(x) = \sum_{n \leq x} \lambda(\mathrm{n})$. Several assertions on classes of summation arithmetic functions, with bounded and unbounded terms of arithmetic functions, which have the property of asymptotic independence, have been proved [7]. In addition, conditions were found [7], under which the limit for summation arithmetic functions is the normal distribution law.

However, only one particular case of summation arithmetic functions having the limit normal distribution law was considered [7], when asymptotically independent summands arithmetic functions have fairly close distribution functions. This case is similar to the normal limit distribution of the sum of independent identically distributed random variables.

In fact, the normal distribution law is much more widespread, for example, for the sum of independent, unequally distributed limited random variables.

However, it was proved in [7] that for bounded arithmetic functions, not independence is satisfied, but only asymptotic independence, more precisely, pairwise asymptotic independence. It may not be sufficient for the summation function with such terms to have a limit normal distribution.

The purpose of this work is to find sufficient conditions under which the summation arithmetic functions with limited terms have a normal distribution law.

## 2. PROPEREITES OF SEQUENCE  RANDOM VARIABLES

Recall some properties of a sequence of random variables $x_i (i = 1,...,n)$.

A sequence of random variables is called stationary in a broad sense, if the conditions [8] are fulfilled:

1. The mathematical expectation $M[x_n]$ does not depend on $n$.          (2.1)

2. $\mathrm{cov}(x_{n+m}, x_n) = \mathrm{M}((\mathrm{x}_{n+m} - M(x_{n+m}))(\mathrm{x}_n - M(x_n)))$ does not depend on $n$.     (2.2)

3. $M[x_n^2] < \infty$ .          (2.3)

We are interested in the asymptotic behavior of a sequence of random variables, i.e. when the value $n \to \infty$.

Поэтому условие (2.1) можно записать в виде:

$$\lim_{n \to \infty} M[x_n] = \mathrm{C},$$
(2.4)

где $C$ - постоянная (возможно равная 0).

Therefore, condition (2.1) can be written as

$$\lim_{n \to \infty} M[x_n] = \mathrm{C},$$
(2.4)

where $C$ is a constant (possibly equal to 0).

A sequence of random variables $x_1, x_2, \ldots$ is called independent if it is executed:

$$P(x_1 \in B_1, \mathrm{x}_2 \in B_2, \ldots) = P(x_1 \in B_1)P(\mathrm{x}_2 \in B_2)\ldots$$
(2.5)

If the random variables $x_k$ and $x_{k+n}$ are asymptotically independent (when the value $n \to \infty$), then based on (2.5) for any values $B_1, B_2 \in B$ is performed:

$$P(x_k \in B_1, \mathrm{x}_{k+n} \in B_2) \to P(x_0 \in B_1)P(\mathrm{x}_0 \in B_2).$$
(2.6)

If the random variables of the sequence $x_i (i = 1, \ldots, n)$ are pairwise asymptotically independent, then the following is true:

$$\lim_{n \to \infty} \mathrm{cov}(x_k, x_{k+n}) = 0$$
(2.7)

Having in mind (2.7), the random variables of the sequence $x_i (i = 1, \ldots, n)$ are pairwise asymptotically independent if the following holds:

$$\lim_{n \to \infty} \mathrm{cov}(x_k, x_{k+n}) = \lim_{n \to \infty} (M[x_k x_{k+n}] - M[x_k]M[x_{k+n}]) = 0.$$
(2.8)

Thus, the enhanced condition (2.2) can be written in the form (2.8).

The last condition (2.3) is written in the form:

$$\lim_{n \to \infty} M[x_n{}^2] < \infty.$$
(2.9)

Now let us talk about the strong mixing property for stationary in a broad sense sequences.

A stationary in a broad sense sequence satisfies the condition of strong mixing, if:

$$\alpha(n) = \sup_{A \in M_1^n, B \in M_{n+k}^{\infty}} |P(AB) - P(A)P(B)| \to 0 \text{ при } n \to \infty, \qquad (2.10)$$

where $\alpha(n)$ is the mixing coefficient, $A \in M_1^n$ - means that the event $A$ belongs to the algebra of events $M_{n+k}^{\infty}$, similarly with respect to the event $B$.

We now apply these properties to the study of summation arithmetic functions.

3. PROPERTIES OF SUMMATION ARITHMETIC FUNCTIONS

As mentioned above, an arbitrary (real) arithmetic function $f$ (or rather, its restriction on $\Omega_n$) in the initial segment of the natural series can be viewed as a sequence of random variables $f_n$ defined on different probability spaces: $f_n : f_n(k) = f(\mathrm{k}), (\mathrm{k} = 1, ..., \mathrm{n})$.

Let us consider property (2.4) for bounded arithmetic functions from these positions. We prove the following assertion.

Assertion 1

Suppose there is an summation function $S(n) = \sum_{k=1}^{n} f(k)$ and the sequence of random variables $f_n$ defined on different probability spaces: $f_n : f_n(k) = f(\mathrm{k}), (\mathrm{k} = 1, ..., \mathrm{n})$. Then, if the arithmetic function $f$ is bounded and there is a limit $\lim_{n \to \infty} \dfrac{S(n)}{n}$ then:

$$\lim_{n \to \infty} M[f_n] = \lim_{n \to \infty} \frac{S(n)}{n} = C, \qquad (3.1)$$

where $C$ is a constant (possibly equal to 0).

Proof

Based on the definition of a random variable, we obtain:

$$M[f_n] = \frac{\sum_{k=1}^{n} f_n(k)}{n} = \frac{\sum_{k=1}^{n} f(k)}{n}.\tag{3.2}$$

Having in mind that $S(n) = \sum_{k=1}^{n} f(k)$, then we obtain based on (3.2):

$$M[f_n] = \frac{S(n)}{n}.\tag{3.3}$$

If the arithmetic function $f$ is limited, then two situations are possible:

1. If there is an asymptotic upper bound $S(n) = o(n)$, then, based on (3.2) and (3.3):

$$\lim_{n\to\infty} M[f_n] = \lim_{n\to\infty} \frac{S(n)}{n} = 0,$$

thus (3.1) with $C = 0$ is satisfied

2. If there are an asymptotic upper bound $S(n) = O(n)$ and $\lim_{n\to\infty} \frac{S(n)}{n}$, then based on (3.2) and (3.3):

$$\lim_{n\to\infty} M[f_n] = \lim_{n\to\infty} \frac{S(n)}{n} = C,$$

thus (3.1) with $C \neq 0$.is satisfied.

Examples of the first case are summation functions: Mertens - $M(n)$, Liouville - $L(n)$ and the number of prime numbers that do not exceed $n$ - $\pi(n)$.

An example of the second case is the summation function of the number of numbers free from squares - $Q(n)$.

Having in mind (3.1), using [1], we get:

$$\lim_{n\to\infty} M[f_n] = \lim_{n\to\infty} \frac{6n/\pi^2 + O(n^{1/2})}{n} = \frac{6}{\pi^2}.$$

However, condition 2) is not necessary and sufficient. For example, for Chebyshev functions containing unbounded terms, it holds $\lim_{n\to\infty}\dfrac{S(n)}{n}=1$.

On the other hand, the limit $\lim_{n\to\infty}\dfrac{S(n)}{n}$ may not exist at all. For example, for an arithmetic function $S(n)=n\mu(n)=\mathrm{O}(n)$, where $\mu(n)$ is the Möbius function.

Now we consider the pairwise asymptotic independence property for arithmetic functions, which corresponds to (2.8).

We assume that an arithmetic function $f$ is pairwise asymptotically independent if for its average values the following relation is satisfied (when the value $n\to\infty$):

$$M_{ij}[f,n]\to M_i[f,n]M_j[f,n],\tag{3.4}$$

where $M_{ij}[f,n]$ is determined by the formula: $\dfrac{\sum\limits_{i=1}^{n}\sum\limits_{j=1(i\neq j)}^{n}f(i)f(j)}{n(n-1)}$, and $M_i[f,n]M_j[f,n]$ is determined by the formula: $\dfrac{(\sum\limits_{k=1}^{n}f(k))^2-\sum\limits_{k=1}^{n}f^2(k)}{n^2}$.

It was proved [7] that for bounded arithmetic functions the pairwise asymptotic independence property (3.4) holds. It is called there simply asymptotic independence for the purpose of contraction.

Let's consider the property of arithmetic functions corresponding to (2.9). We prove the following assertion.

Assertion 2

If we have a bounded arithmetic function $f$, then the following relation holds:

$$\lim_{n\to\infty}M[f^2,n]<\infty,\tag{3.5}$$

where $M[f^2,n]$ - the average value of the arithmetic function $f^2$.

Proof

Having in mind that the arithmetic function $f$ is bounded, then $|f| \leq L$, where $L$ is a constant, therefore it is executed:

$$M[f^2, n] = \frac{\sum_{k=1}^{n} f^2(k)}{n} \leq L^2 < \infty,$$

which corresponds to (3.5).

The sequence members $f_n$ are in different probability spaces, and the stationarity property in a broad sense is introduced for a sequence of random variables whose members are in the same probability space. Therefore, we call this property for a sequence $f_n$ pseudo stationarity in a broad sense.

Now we consider the strong mixing property for summation arithmetic functions.

For example, let us consider two summation arithmetic functions: $\pi(n), Q(n)$, for which the corresponding random sequences satisfy the conditions of pseudo-stationarity in a broad sense.

Based on (2.10), we find the strong mixing coefficient for a random sequence corresponding to an arithmetic function of the number of primes. We take the event that the number $k$ is prime as the event $A$. We take the event $B$ that the number $n + k$ is composite, i.e. $x_{n+k} = 0$. Then $P(x_k = 1) = 1/\ln k + o(1/\ln k)$, $P(x_{n+k} = 0) = 1 - 1/\ln(n+k) + o(1/\ln(n+k))$ and $P(A)P(B) = (1 + 1/\ln k + o(1/\ln k))(1 - 1/\ln(n+k) + o(1/\ln(n+k)))$. On the other hand $P(AB) = P(A)P(B/A) = 1/\ln k + o(1/\ln k)$, since $P(B/A) = 1$, if the number $n + k$ is even. The last assumption is made to obtain the upper limit of the difference. Therefore, the coefficient of strong mixing in this case is equal to:

$$\alpha(n) = \sup |P(AB) - P(A)P(B)| = (1/\ln k + o(1/\ln k))(1/\ln(n+k) + o(1/\ln(n+k))) .(3.6)$$

Having in mind (3.6) $\alpha(n) \to 0$ when $n \to \infty$, therefore, the condition of strong mixing is satisfied in this case.

Now check the condition $\sum_{n=1}^{\infty} \alpha(n) < \infty$ for the number of primes:

$$\sum_{n=1}^{\infty} \alpha(n) = [1/\ln k + o(1/\ln k)] \sum_{n=1}^{\infty} 1/\ln(n+k) = \infty.$$  (3.7)

Based on (3.7), this condition is not satisfied for the arithmetic function of the number of primes.

Now we find the strong mixing coefficient for a sequence of random variables corresponding to the arithmetic function of the number of natural numbers free from squares.

Having in mind [1], the following asymptotic formula is valid:

$$Q(n) = 6n/\pi^2 + O(n^{1/2}).$$  (3.8)

We take the event that the number $k$ is square-free as an event $A$, i.e. $x_k = 1$. We take the event that the number is not free squares as an event $B$, i.e. $x_{n+k} = 0$. Then, based on (4.3): $P(x_k = 1) = 6/\pi^2 + O(1/k^{1/2})$, $P(x_{n+k} = 0) = 1 - 6/\pi^2 + O(1/(n+k)^{1/2})$. The product of the probabilities is $P(A)P(B) = [6/\pi^2 + O(1/k^{1/2})][1 - 6/\pi^2 + O(1/(n+k)^{1/2})]$ in this case. Then the value $P(AB)$ is: $P(A)P(B/A) = [6/\pi^2 + O(1/k^{1/2})][1 - 6/\pi^2 + O(1/(n+k)^{1/2})]$. The last equality is explained by the independence of the event $B$ from the event $A$.. Therefore, the condition of strong mixing is fulfilled in this case:

$$\alpha(n) = \sup |P(AB) - P(A)P(B)| = 0.$$  (3.9)

Based on (3.9), the condition $\sum_{n=1}^{\infty} \alpha(n) = 0 < \infty$ is satisfied for a sequence of random variables corresponding to a given arithmetic function.

Since the condition of strong mixing is introduced for one probabilistic space, then the fulfillment of this condition for the sequence $x_n$ is called the condition of pseudo strong mixing of a given sequence.

4.  SUFFICIENT CONDITION TO HAVE A LIMITING NORMAL DISTRIBUTION LAW FOR A SUMMATION ARITHMETIC FUNCTION WITH BOUNDED TERMS

The following theorem was proved in [9].

Let we have a stationary in a broad sense sequence of random variables $x_n (n = 1, 2, ...)$ satisfies the property of strong mixing with $\sum_{n=1}^{\infty} \alpha(n) < \infty$. Suppose, moreover, that the random variables for given sequence are bounded with probability 1 - $|x_n| = c_o < \infty$. Then, if $D(x_n) \neq 0$, then the random variable $S_n = \sum_{k=1}^{n} x_k$ has a limit normal distribution.

We cannot directly apply this theorem, since, according to the conditions of this theorem, a sequence of random variables $x_n (n = 1, 2, ...)$ is defined in one probability space and in our case, random values $f_n(k) = f(k), (k = 1, ..., n)$ are defined for different values $n$ in different probability spaces.

Based on Lemma 3, p. 123 [10], if a sequence of random variables $f_n : f_n(k) = f(k)$ defined in different probability spaces converges in distribution $P(f_n < x) \rightarrow P(f < x)$ when the value $n \rightarrow \infty$, then a sequence of random variables $g_n$ (in one probability space) can be constructed that converges in distribution $P(g_n < x) \rightarrow P(g < x)$ when the value $n \rightarrow \infty$ with $P(f_n < x) = P(g_n < x)$, $P(f < x) = P(g < x)$.

We will perform the transformation corresponding to the indicated lemma for random variables that are in different probability spaces. The distribution functions of random variables and their expected values are preserved under this transformation. Thus, with this transformation, the pseudo asymptotic independence of random variables is preserved, when the value $n \rightarrow \infty$:

$$M[x_k x_{k+n}] - M[x_k] M[x_{k+n}] \rightarrow 0. \qquad (4.1)$$

We can say that the independence of random variables $x_k$ and $x_{k+n}$ is performed when $M[x_k x_{k+n}] = M[x_k] M[x_{k+n}]$, except in the trivial cases, when at least one of the mathematical expectations is 0.

As mentioned above, another formulation of asymptotic independence, which also includes trivial cases, when $n \rightarrow \infty$ a value, is:

$$P(x_k \in B_1, x_{k+n} \in B_2) \rightarrow P(x_k \in B_1) P(x_{k+n} \in B_2). \qquad (4.2)$$

Therefore, the formulations of asymptotic independence (4.1) and (4.2) are equivalent among themselves everywhere except in the trivial cases.

Based on the definition of a random variable, it follows that for sets $B_1, B_2$ from the sigma of an algebra of Borel sets on a line, there are: $A = x_k^{-1}(B_1), B = x_{k+n}^{-1}(B_2)$, where $B_1 \in M_1^n, B_2 \in M_{n+k}^\infty$.

Therefore, from performing (4.2) it follows the fulfillment of the strong mixing property (2.10) when the value $n \to \infty$:

$$\alpha(n) = \sup\nolimits_{A \in M_1^n, B \in M_{n+k}^\infty} \mid P(AB) - P(A)P(B) \mid \to 0.$$

Having in mind the above, we can formulate the following assertion.

Assertion 3

Suppose there is a summation arithmetic function $S(n) = \sum_{k=1}^n f(k)$, where $\lim_{n \to \infty} \dfrac{S(n)}{n} \neq 0$, and a sequence of random variables $f_n : f_n(k) = f(k)$ defined in different probability spaces (respectively, with distribution functions $P(f_n < x)$, which converges in distribution $P(f_n < x) \to P(f < x)$ when the value $n \to \infty$, is pseudo stationary in a broad sense, satisfies pseudo strong mixing conditions with $\sum_{n=1}^\infty \alpha(n) < \infty$ and is bounded. Then, if $D[f_n] \neq 0$, then the summa arithmetic function $S(n) = \sum_{k=1}^n f(k)$ has a limit normal distribution law.

Proof

So, random variables corresponding to an arithmetic function $f$ — $f_n : f_n(k) = f(\text{k}), (\text{k} = 1, ..., \text{n})$ for different values $n$, are defined in different probability spaces.

We use the transformation corresponding to Lemma 3, p. 123 [10] for random variables that are in different probability spaces.

In the particular case, as shown above, this transformation is performed for a sequence of random variables $f_n : f_n(k) = f(k), (k = 1, ..., n)$ corresponding to an arithmetic function $f$ located in different probability spaces.

Let us suppose that after this transformation a sequence of random variables $g_n$ is in the same probability space. Then, based on indicated lemma, the distribution of random variables $g_n$, respectively, coincides with the distribution of random variables $f_n$.

Accordingly, the random variable $G_n = \sum_{k=1}^{n} g_k$ has the same distribution function as the random variable $S_n : S_n(k) = S(k), (k = 1, ..., n)$ - $P(G_n < x) = P(S_n < x)$ and the limiting distribution functions also coincide - $P(G < x) = P(S < x)$.

This transformation also preserves the pseudo-asymptotic independence of random variables (4.1), which is performed for a sequence of bounded random variables $f_n$ corresponding to a bounded arithmetic function $f$ [7], and, therefore, the pairwise asymptotic independence of random variables of the sequence $g_n$ in one probability space is performed when the value $n \to \infty$:

$$M[g_k \ g_{k+n}] - M[g_k]M[g_{k+n}] \to 0. \tag{4.3}$$

As was shown above, pairwise asymptotic independence of random variables (4.3) is equivalent to asymptotic independence (4.2) when the condition is satisfied (no trivial cases), and therefore, this transformation preserves the strong mixing factor of these random variables when the value $n \to \infty$:

$$\alpha(n) = \sup_{A \in M_1^n, B \in M_{n+k}^\infty} | P(AB) - P(A)P(B) | \to 0. \tag{4.4}$$

According to the condition, the sequence of random variables $f_n$ is stationary in a broad sense and satisfies the condition of strong mixing, therefore the sequence of random variables $g_n$ is also stationary in a broad sense and satisfies the condition of strong mixing, since the distributions of random variables and, accordingly, the characteristics coincide.

Based on $D[f_n] = D[g_n] \neq 0$ and theorem the random variable $G$ has a normal distribution.

Having in mind that the distribution of random variables are the same: $P(G < x) = P(S < x)$, a random variable $S : S_n \to S$ (by distribution when the value $n \to \infty$) also has a normal distribution. Based on $S_n : S_n(k) = S(k), (k = 1, ..., n)$ the summation function $S(n) = \sum_{k=1}^{n} f(k)$ has a limit normal distribution.

Let us verify the fulfillment of the conditions of this statement for summation arithmetic with bounded terms.

The condition $D[f_n] \neq 0$ corresponds to such condition $D[f] \neq 0$ for the arithmetic function $f$. We consider only such arithmetic functions, because if the arithmetic function $D[f] = 0$ deviates from its mean value, it is 0 and such arithmetic functions are of no interest to us.

An asymptotic estimate is performed: $M(n) = o(n), \text{L}(n) = \text{o}(n), \pi(n) = \text{o}(n)$ for the considered summation arithmetic functions with bounded terms, therefore $\lim_{n \to \infty} S(n) / n = 0$ for the indicated functions and they do not satisfy the sufficient condition of assertion 3.

На основании [1] для суммарной арифметической функции количество натуральных чисел, свободных от квадратов и не превосходящих значение $n$ справедливо: $\lim_{n \to \infty} Q(n) / n = 6 / \pi^2 \neq 0$, поэтому данная суммарная функция удовлетворяет данному достаточному условию утверждения 3.

Based on [1] for the summation arithmetic function - the number of natural numbers free from squares and not exceeding the value $n$ is true. Therefore this summation function satisfies this sufficient condition of assertion 3.

Previously it was shown that the pseudo strong mixing property is fulfilled for the summation arithmetic function $Q(n)$, therefore, based on assertion 3, the limit distribution for this function is the normal distribution. This is consistent with the QEIS A158819 schedule.

The condition $\lim_{n \to \infty} S(n) / n \neq 0$, as mentioned earlier, is sufficient, but not necessary, so this case needs additional investigation. For example, summation functions $\sum_{k=1}^{n} \frac{\mu(k)}{k}, \sum_{k=1}^{n} \frac{\lambda(k)}{k}$, where $\mu(k), \lambda(k)$ are, respectively, the Möbius and Liouville functions at the point $k$, have a limit normal distribution law, despite the fact that $\lim_{n\to\infty} \frac{1}{n}\sum_{k=1}^{n}\frac{\mu(k)}{k} = \lim_{n\to\infty}\frac{1}{n}\sum_{k=1}^{n}\frac{\lambda(k)}{k} = 0$. This case is considered in more detail in [7].

### 3. CONCLUSION AND SUGGESTIONS FOR FURTHER WORK

The next article will continue to study the behavior of some sums.

### 4. ACKNOWLEDGEMENTS

Thanks to everyone who has contributed to the discussion of this paper. I am grateful to everyone who expressed their suggestions and comments in the course of this work.